\input amstex
\documentstyle{amsppt}
\magnification1100
\topmatter
\title Braided near-group categories \endtitle
\author Jacob Siehler \endauthor

\address Department of Mathematics, Virginia Tech, Blacksburg VA 24061-0123 USA\endaddress
\email jsiehler\@math.vt.edu \endemail


\abstract
A near-group category is an additively semisimple category with a product
such that all but one of the simple objects is invertible. We classify braided 
structures on near-group categories, and give explicit numerical formulas for
their associativity and commutativity morphisms. \endabstract

\thanks Partially supported by the National Science Foundation \endthanks

\endtopmatter

\document 

\def\hom{\hbox{hom}}

\def\inv{^{-1}}
\def\C{ \Cal C}
\def\invariant{(\delta_1,\ldots,\delta_n,\epsilon)}

\head 1. Results \endhead
We consider near-group categories; that is, semisimple monoidal categories, 
whose set of simples consists of a finite group $G$ of invertible elements together with
one noninvertible element $m$. We restrict to the case where $m$ is not a summand in $mm$,
and we will say that $G$ is the underlying group for such a category.  Our categories
also have a ground ring $R$, and we assume throughout that $R$ is an integral
domain.

Tambara and Yamagami \cite{TY} have classified the possible monoidal structures on such
categories.  Their result is:
\hfuzz2.3pt
\proclaim {1.1 Theorem (Tambara-Yamagami)}
Monoidal near-group categories correspond to pairs $(\chi,\tau)$ where
$\chi$ is a nondegenerate, symmetric $R$-valued  bicharacter of $G$,
and $\tau$ is a square root of $1/|G|$
\endproclaim
\hfuzz0.1pt
Our main result is the following theorem, which gives necessary and
sufficient conditions under which these categories admit a balanced (tortile)
braided structure, and parametrizes the distinct braidings.

\par
\proclaim {1.2 Theorem}
Suppose $G$ is a finite abelian group, $\C$ is a monoidal near-group category
with underlying group $G$, and $R$ contains roots of unity of order $8|G|$.
Then

(1) $\C$ admits a braiding if and only if $G$ is an elementary abelian 2-group; that is,
every element has order~2. \par
(2) The nonequivalent braidings on $\C$ are in one-to-one correspondence with $(n+1)$-tuples $(\delta_1,\ldots,\delta_n,\epsilon)$, where $\epsilon=\pm1$ and $\delta_i=\pm1$ for all $i$,
and $n$ is the rank of $G$.\par
(3) Each braiding of $\C$ has exactly two choices of twist morphisms
compatible with it.\par
\endproclaim
Explicit computation of the commutativities further establishes

\proclaim {1.3 Corollary} In any braiding of $\C$,\par
(1) the group subcategory generated by invertibles is symmetric commutative.\par
(2) a basis can be chosen so that every commuting isomorphism is multiplication by 
an $8|G|$-th root of unity.

\endproclaim
In section 5 we give a simple example to show that that the order
condition of 1.3(2) is sharp.  Section 6 contains some
additional comments about the conditions on the ground ring $R$.

\subhead 1.4 Remarks\endsubhead
\roster
\item It is easily seen that if the rank of $G$ is odd then the category
does not have an
integer-valued dimension function. In particular it is not a category of
representations.
\cite{TY} shows that some even-rank cases are not representation
categories, for more
subtle reasons.
\item According to \cite{TY}, any abelian group can underlie an associative
near-group category. Similarly \cite{Q} shows any abelian group underlies a braided
group-category. The theorem and corollary show that the structures that extend are quite
rare: few associative near-group categories support a braided structure, and few braided
group-categories can be embedded in a braided near-group category.
\item
This rarity may be related to duality. The self-dual elements of a 
group-category correspond exactly to the elements of order two. Our definition of
near-group category requires the object $m$ to be self-dual, and the conclusion is that this
forces the invertible elements to be self-dual too. Generalizing the definition to
allow two mutually dual non-invertible simple objects might enlarge the number of groups allowed.
\endroster       
         
\head 2.  The invariants \endhead

\subhead 2.1 Notation \endsubhead
In a near-group category there are four different kinds of commutes to
consider, and we use notation that recognizes
their different roles.  For $g,h\in G$, the commuting isomorphisms are:
$$\halign{\indent#&#&#\hfil\cr
$gh \to hg$&\quad is multiplication by &$\sigma_0(g,h)$ \cr
$gm \to mg$&\quad is multiplication by &$\sigma_1(g)$ \cr
$mg \to gm$&\quad is multiplication by &$\sigma_2(g)$ \cr
$mm \to mm$&\quad is multiplication by &$\sigma_3(g)$ on the $g$ summand\cr}$$

\subhead 2.2 Extracting the invariants \endsubhead
We first want to clarify the result of \cite{TY} cited in 1.1. Saying that $\chi$ is a bicharacter of $G$ means that $\chi(ab,c)=\chi(a,c)\chi(b,c)$  and $\chi(a,bc)=\chi(a,b)\chi(a,c)$;  this is just a bilinear form written multiplicatively.
The technique of [TY] is to establish that for each pair $(\chi,\tau)$
there is a choice of basis, called a {\it normal basis},
so that the associativities in the category have the following form:
for $a,b,c\in G$,
$$\eqalign{
\alpha_{a,b,c}&=1\cr
\alpha_{a,b,m}=\alpha_{m,a,b}&=1\cr
\alpha_{a,m,b}&=\chi(a,b)\cr
\alpha_{a,m,m}=\alpha_{m,m,a}&=\oplus_b 1_b\cr
\alpha_{m,a,m}&=\oplus_b \chi(a,b)\cr
\alpha_{m,m,m}&=\left(\tau\chi(a,b)\inv\right)_{a,b}\cr
}$$
and moreover,
$$\alpha_{m,m,m}\inv=\left(\tau\chi(a,b)\right)_{a,b}$$

Now, to obtain the $\invariant$ describing the commutativity, fix generators $g_1$ through $g_n$ for $G$ and choose a normal basis for the parameters $\chi$ and $\tau$ describing the associative structure.  Then the $\delta_i$ and $\epsilon$ parameters are computed by taking
$$\delta_i={\sigma_1(g_i) \over \sqrt{\chi(g_i,g_i)}}$$
and
$$\epsilon={\sigma_3(1) \over \sqrt{\tau\sum_{g\in G}\sigma_1(g)}}$$
It is not obvious that these invariants are $\pm1$-valued;
this will follow from our examination of the hexagon equations below.

\subhead 2.3 Realizing the parameters \endsubhead
[TY] gives the construction to realize the parameters $\chi$ and $\tau$ describing the associative structure.  We assume we have that construction
in hand and fix a normal basis for it.  Realizing the braiding parameters
$(\delta_1,\ldots,\delta_n,\epsilon)$ means explicitly giving the values of the 
$\sigma$ functions.\par
1.  Commuting group elements. Set $$\sigma_0(g,h)=\chi(g,h)$$\par\medskip
2.  Commuting a group element with the noninvertible.
First, for the generators $g_i$, set
$$\sigma_1(g_i)=\sigma_2(g_i)=\delta_i\sqrt{\chi(g_i,g_i)}$$
Any other element $g$ in the group can be written uniquely as a product
of generators with each generator appearing at most once (generators
have order 2).  Say $g=h_1h_2\cdots h_k$ is the expression for $g$ as
a product of generators.  Set
$$\sigma_1(g)=\sigma_2(g)=\prod_{i=1}^{k}
 \left(\sigma_1(h_i)\prod_{i<j\le k}\chi(h_i,h_j)\right)$$
\par\medskip

3.  Commuting the noninvertible with itself.  This also goes in two steps.  First, set
$$\sigma_3(1)=\epsilon\sqrt{\tau\sum_{g\in G}\sigma_1(g)}$$
Then, for any other $g\in G$, set
$$\sigma_3(g)=\sigma_3(1)\sigma_1(g)\chi(g,g)$$
\par
Comparing these definitions to the procedure above for extracting invariants, it is clear that the $\sigma$'s we have defined do realize the parameters$(\delta_1,\ldots,\delta_n,\epsilon)$. \par

\subhead 2.4 Note \endsubhead
We see from the formulas that the commutativity inside the group category is entirely determined by the associativity parameters; all the ``new" information concerns the behavior of the noninvertible.

\head 3. Proofs \endhead

\subhead 3.1 The hexagon restraints \endsubhead 
The results of the theorem come directly from the hexagon axioms.  If we let $a,b,c$ stand for group elements, there are eight types of three-term products to consider: $abc$, $abm$, $amb$, $mab$, $mma$, $mam$, $amm$, and $mmm$.  Using a normal basis, we can express associativities in terms of $\chi$ and $\tau$, and write out the standard hexagon equations corresponding to the eight three-term products:
$$\eqalign{
\sigma_0(a,b)\sigma_0(a,c)&=\sigma_0(a,bc)\cr
\sigma_0(a,b)&=\chi(b,a)\cr
\sigma_0(a,b)\sigma_1(a)&=\sigma_1(a)\chi(a,b)\cr
\sigma_2(ab)&=\sigma_2(b)\chi(a,b)\sigma_2(a)\cr
\sigma_2(a)\sigma_3(ba\inv)&=\chi(a,b)\sigma_3(b)\cr
\sigma_3(a\inv b)\sigma_2(a)&=\sigma_3(b)\chi(a,b)\cr
\sigma_0(a,a\inv b)&=\sigma_1(a)\chi(a,b)\sigma_1(a)\cr
\sigma_3(a)\tau\chi(a,b)\inv\sigma_3(b)
    &=\sum_{c\in G}\tau\chi(a,c)\inv\tau\sigma_2(c)\chi(c,b)\inv\cr
}
$$

The hexagon for $mmm$ actually gives the matrix equation $AS_2A=S_3AS_3$, where 
$A=\alpha_{m,m,m}=(\chi(a,b))_{a,b}$, $S_2=\oplus_a \sigma_2(a)$ and $S_3=\oplus_a \sigma_3(a)$.
The final equation in the list above corresponds to the $(a,b)$-th entry of the matrix equation.\par

There are some obvious redundancies among the hexagon equations; eliminating these and simplifying the remaining ones yields the following reduced set of equations, equivalent to the original eight:

$$\eqalignno{
\sigma_0(a,b)&=\chi(a,b)&(1)\cr
\sigma_2(ab)&=\sigma_2(a)\sigma_2(b)\chi(a,b)&(2)\cr
\sigma_2(a)\sigma_3(ba\inv)&=\chi(a,b)\sigma_3(b)&(3)\cr
\sigma_1(a)^2&=\chi(a\inv,a)&(4)\cr
\tau\sum_c\chi(ab,c)\inv\sigma_2(c)&=\chi(a,b)\inv\sigma_3(a)\sigma_3(b)&(5)\cr}$$

Computing the inverse hexagons amounts to inverting associativities and reversing commutativities, producing:

$$\eqalignno{
\sigma_0(b,a)&=\chi(a,b)\inv&(6)\cr
\sigma_1(ab)&=\sigma_1(a)\sigma_1(b)\chi(a,b)\inv&(7)\cr
\sigma_1(a)\sigma_3(ba\inv)&=\chi(a,b)\inv\sigma_3(b)&(8)\cr
\sigma_2(a)^2&=\chi(a,a)&(9)\cr
\tau\sum_c\chi(ab,c)\sigma_1(c)&=\chi(a,b)\sigma_3(a)\sigma_3(b)&(10)\cr}$$

\subhead 3.2 Proof of 1.2(1) \endsubhead
Equations (1) and (6) together with the fact that $\chi$ is symmetric imply that $\chi(a,b)=\chi(a,b)\inv$ for every $a,b\in G$; that is, $\chi$ is $\pm1$-valued.
But $\chi$ is also nondegenerate.  If $|a|=2n+1$, then $\chi(a,b)^{2n+1}=\chi(1,b)=1$
for every $b$.  So $\chi(a,b)=1$ for every $b$, contradicting nondegeneracy.  If $|a|=2n$ for
$n>1$ then $a^2\ne1$ but $\chi(a^2,b)=1$ for every $b$, also contradicting nondegeneracy.  This shows that every element in $G$ has order 2, as claimed.\par

\subhead 3.3 Note \endsubhead
Since group elements are self-inverse, we will now drop the inverse signs from group
elements whenever they appear.  Likewise, we replace $\chi\inv$ with $\chi$ wherever it appears.\par

\subhead 3.4 Proof of Corollary 1.3(1) \endsubhead
For invertibles $g$ and $h$, 
$$\sigma_0(g,h)\sigma_0(h,g)=\chi(g,h)\chi(h,g)=\chi(g,h)^2=1.$$

\subhead 3.5 Proof of 1.2(2) \endsubhead
\par
{ 1. The invariant $\invariant$ is well-defined.}
We have to show that our definition is independent of the choice of normal basis that was made.
For any group element $a$, let $a_l$ and $a_r$ be the basis for $\hom(m,am)$ and $\hom(m,ma)$, respectively, in some normal basis; let $a_l'$ and $a_r'$ play the same role in some other normal basis.  With respect to any  normal basis (by definition), the association $(ma)m\to m(am)$ is represented by the identity on the 1 summand.  That implies that there is some $k$ in the ground ring so that $a_l'=ka_l$ and $a_r'=ka_r$; it follows that the commuting isomorphism $\sigma_1(a):am\to ma$ is represented the same with respect to either normal basis, so setting $\delta_i=\sigma_1(g_i)/\sqrt{\chi(g_i,g_i)}$ is well-defined.\par
It is clear that $\sigma_3(1)$ doesn't depend on the choice of basis of $\hom(1,mm)$; together 
with the preceding argument, this is enough to establish that the invariant $\epsilon$ is
well-defined.\par
\medskip
{ 2. Giving $\invariant$ determines the commutativity.}
We have already noted that $\sigma_0$ is totally determined by the associative structure.
Equations (3) and (8) together imply that $\sigma_1\equiv\sigma_2$. Equations (2) and (7) say that
$\sigma_1$ and $\sigma_2$ are determined by their values on the generators $g_1,\ldots,g_n$ for $G$; finally, (4) and (9) say that $\sigma_1(g_i)=\sigma_2(g_i)$ is determined, up to 
a sign, by the value of $\chi(g_i,g_i)$.  Therefore, giving the $\delta_i$ to specify those signs
is enough to determine $\sigma_1$ and $\sigma_2$.\par

Equation (3) says that $\sigma_3$ is determined by its value on any one group element.  For instance, specifying $\sigma_3(1)$ is enough to describe $\sigma_3$ entirely.  But evaluating
(5) for the special case $a=b=1$, we get
$$\sigma_3(1)^2=\tau\sum_{c\in G}\sigma_1(c)$$
So $\sigma_3(1)$ is in fact determined up to a sign, and giving $\epsilon$ to specify that sign is enough to determine $\sigma_3$, completing the description of all the commutativities.\par
\medskip
{ 3. Every $\invariant$ can be realized.}
Define the $\sigma$ functions according to the construction given in 2.3.  We have to show that
all the hexagon equations are satisfied. Clearly, equations (1) and (6) are satisfied.  \par
Next, we will verify (2) and (7); since the construction defines $\sigma_1$ and $\sigma_2$ to
be equal, and since we can drop inverses, it suffices to check that (2) is satisfied.  
For a generator $g_i$,
$$\sigma_2(g_i)^2=\left(\delta_i\sqrt{\chi(g_i,g_i)}\right)^2=\chi(g_i,g_i)$$
so (2) is satisfied on generators.  For any element $g=h_1h_2\cdots h_k$ expressed
as a product of generators,
$$\eqalign{
\sigma_2(g)^2&=\left(\prod_{i=1}^{k}
 \left(\sigma_2(h_i)\prod_{i<j\le k}\chi(h_i,h_j)\right)\right)^2\cr
&=\prod_{i=1}^{k}\sigma_2(h_i)^2\quad\hbox{since $\chi^2=1$}\cr
&=\prod_{i=1}^{k}\chi(h_i,h_i)\cr}
$$
Working from the other end,
$$\eqalign{
\chi(g,g)&=\chi(h_1h_2\cdots h_k,h_1h_2\cdots h_k)\cr
&=\prod_{1\le i,j\le k}\chi(h_i,h_j)\cr
&=\prod_{i=1}^{k}\chi(h_i,h_i)\quad
\hbox {since $\chi(h_i,h_j)\chi(h_j,h_i)=1$ when $i\ne j$}\cr}$$
This establishes that $\sigma_2(g)^2=\chi(g,g)$ for every $g\in G$, so equations $(2)$ and $(7)$
are satisfied.\par\medskip
To show that (3) and (8) are satisfied, it suffices to check one or the other;
we show that (8) holds (remember that we can drop inverses):
$$\eqalign{
\sigma_1(a)\sigma_3(ab)&=\sigma_1(a)\sigma_3(1)\sigma_1(ab)\chi(ab,ab)\cr
&=\sigma_1(a)\sigma_3(1)\sigma_1(a)\sigma_1(b)\chi(a,b)\chi(ab,ab)\cr
&=\sigma_1(a)^2\sigma_3(1)\sigma_1(b)\chi(a,b)\chi(a,a)\chi(b,b)\cr
&=\chi(a,a)^2\chi(a,b)\sigma_3(1)\sigma_1(b)\chi(b,b)\cr
&=\chi(a,b)\sigma_3(b)\cr}$$
Finally we check that (10) is satisfied; (5) is identical.  This is
another computation directly from the definitions:
$$\eqalign{
\chi(a,b)\sigma_3(a)\sigma_3(b)
&=\chi(a,b)\chi(a,a)\chi(b,b)\sigma_1(a)\sigma_1(b)\sigma_3(1)^2\cr
&=\chi(a,b)\chi(a,a)\chi(b,b)\sigma_1(a)\sigma_1(b)\tau\sum_c\sigma_1(c)\cr
&=\tau\sum_c\chi(a,a)\chi(b,b)\sigma_1(ab)\sigma_1(c)\cr
&=\tau\sum_c\chi(a,a)\chi(b,b)\chi(ab,c)\sigma_1(abc)\cr
&=\tau\sum_c\chi(ab,abc)\sigma_1(abc)=\tau\sum_c\chi(ab,c)\sigma_1(c)\cr
}$$
which shows that (10) and (5) are satisfied and completes the verification
that the construction gives a commutativity satisfying all the hexagon
equations.\par

\subhead 3.7 Proof of 1.2(3) \endsubhead
The balance axiom requires that we produce automorphisms $\theta_s$ for each simple object $s$
satisfying $\theta_{rs}=\theta_r\theta_s\sigma(r,s)\sigma(s,r)$.  Here, that means
that for all $g,h\in G$, we must have
$$\eqalign {
\theta_{gh}&=\theta_g\theta_h\sigma_0(h,g)\sigma_0(g,h)\cr
\theta_m&=\theta_g\theta_m\sigma_1(g)\sigma_2(g)\cr
\hbox{and}\quad
\theta_g&=\theta_m^2\sigma_3(g)^2 }$$
We know that commutativity in the group subcategory is symmetric, and with respect to a normal basis we know that $\sigma_1=\sigma_2$.  So we can give a slightly simplified statement of the requirements:
$$\eqalignno {
\theta_{gh}&=\theta_g\theta_h&(11)\cr
\theta_g&=\sigma_1(g)^2&(12)\cr
\theta_g&=\theta_m^2\sigma_3(g)^2&(13) }$$
Strictly speaking, equation (12) should be $\theta_g=\sigma_1(g)^{-2}$,
but $\sigma_1(g)^2=\chi(g,g)=\pm1$, so dropping the inverse is harmless.
To show that the category is balanced we define
$$\eqalign{
\theta_g&:=\sigma_1(g)^2\cr
\hbox{and}\quad \theta_m&:=\pm1/\sigma_3(1)\cr}$$
and check that these satisfy the axioms.  Equation (12) is obviously satisfied.  We
check (11):
$$\theta_{gh}=\sigma_1(gh)^2=\left(\sigma_1(g)\sigma_1(h)\chi(g,h)\right)^2
=\sigma_1(g)^2\sigma_1(h)^2=\theta_g\theta_h$$
and (13):
$$\eqalign{
\theta_m^2\sigma_3(g)^2&={1\over\sigma_3(1)^2}(\sigma_3(1)\sigma_1(g)\chi(g,g))^2\cr
&=\sigma_1(g)^2=\theta_g\cr}$$
This shows that there are two choices of twist morphisms which balance the commutativity;
conversely, it is obvious that we had to define the twist morphisms as we did, so these
are the only two possibilities, concluding the proof of Theorem 1.
\subhead 3.8 Proof of Corollary 1.3(2) \endsubhead
Assume that the $\sigma$ functions are defined according to the construction in 2.3.  We want to show that they take values in the roots of unity.  First, $\sigma_0(g,h)=\chi(g,h)$ so $\sigma_0$ is $\pm1$-valued.
Also, $\sigma_1(g)^2=\sigma_2(g)^2=\chi(g,g)$ so $\sigma_1$ and $\sigma_2$ take
values that are at most 4th roots of unity.  \par
Now, recall from 3.1 that the hexagon for $mmm$ can be expressed as the matrix equation $AS_2A=S_3AS_3$.
Pass to determinants; $A^2=I$, so $\det A=\pm1$, and
$$\eqalign{
\det A\det S_1&=\det S_3^2\cr
\pm\prod_g\sigma_1(g)&=\prod_g\sigma_3(g)^2\cr
\pm\prod_g\sigma_1(g)&=\prod_g\sigma_3(1)^2\sigma_1(g)^2\chi(g,g)^2\cr
\sigma_3(1)^{2|G|}&=\pm\prod_g \sigma_1(g)\inv\cr
}$$
so $\sigma_3(1)^{8|G|}=1$.  For any other
$g\in G$, $\sigma_3(g)=\sigma_3(1)\sigma_1(g)\chi(g,g)$, so $\sigma_3(g)^{8|G|}=1$, completing the proof.
\medskip
\head 5. An example \endhead
If we let $G=\langle g\rangle$ be the cyclic 2-group, we obtain examples showing that the bounds
of Corollary 1.3(2) are best possible without additional hypotheses.  There is only one nondegenerate bicharacter $\chi$ on $G$; monoidal structures are given by $\tau=\pm(1/\sqrt2)$, and braidings are given by pairs $(\delta,\epsilon)$ of $\pm1$'s.  If we leave the parameters generic, we can compute the interesting commutativities ($I$ stands for a primitive 4th root of 1):
$$\eqalign{
\sigma_1(1)&=\sigma_2(1)=1\cr
\sigma_1(g)&=\sigma_2(g)=\delta I\cr
\sigma_3(1)&=\epsilon\sqrt{\tau(1+\delta I)}\cr
\sigma_3(g)&=\epsilon\delta I\sqrt{\tau(1+\delta I)}\cr}$$
In this case, regardless of the parameters, it is easily verified that $\sigma_3(1)$ and $\sigma_3(g)$ are primitive 16th roots of unity.
\par\medskip
\head 6.  The ground ring \endhead
Theorem 1.2 says that if $R$ has enough roots of unity, then every
$\invariant$ can be realized as a braiding over $R$.  In general,
if $R$ does not have enough roots of unity, some of these will
be impossible to realize.  The example from section 5 should be
sufficiently cautionary: there, none of the parameter sets can be
realized unless the exact hypotheses of theorem 1.2 are satisfied.
On the other hand, the explicit formulas we have for the associativies
and commutativities show that a lack of roots of unity is
essentially the only obstruction to realizing the braidings.\par
  We note that near-group categories are never defined over the integers;
this is ruled out just by the monoidal structure, since $\tau$ is never
integral, yet it necessarily appears in the associativity of the
noninvertible. We also point out that the bicharacter $\chi$ has an important
effect on what roots of unity are required.  If $\chi(g,g)=-1$
for any $g\in G$, then at least 4th roots of unity are unavoidable in realizing
the braidings (this is because $\sigma_1(g)^2=\chi(g,g)$). Conversely, 
if $\chi(g,g)=+1$ for every $g$, then we can refine the proof of 1.3(2) 
to show that we need only $4|G|$-th roots of unity to realize all the braidings.

\enddocument